\pdfoutput=1
\documentclass[12pt,reqno]{amsart}

\usepackage[utf8]{inputenc}
\usepackage[T1]{fontenc}
\usepackage[charter]{mathdesign}
\usepackage[draft=false,kerning=true]{microtype}

\usepackage{fullpage}

\usepackage{amsthm, amsmath, amsfonts}
\usepackage{mathtools}
\usepackage{enumerate}

\usepackage{todonotes}

\usepackage{graphicx}
\usepackage{tikz}
\usetikzlibrary{arrows,backgrounds,calc}
\usepackage{pgfplots}
\pgfplotsset{compat=1.12}

\usepackage{hyperref}

\theoremstyle{remark}

\theoremstyle{definition}

\newtheorem*{example*}{Example}

\renewcommand{\MR}[1]{}

\newmuskip\pFqskip
\mathchardef\pFcomma=\mathcode`, % keep a copy of the comma
\newcommand*\pFq[5]{
  \begingroup
  \begingroup\lccode`~=`,
    \lowercase{\endgroup\def~}{\pFcomma\mkern\pFqskip}%
  \mathcode`,=\string"8000
  {}_{#1}F_{#2}\biggl(\genfrac..{0pt}{}{#3}{#4}\bigg|#5\biggr)%
  \endgroup
}

\DeclareSymbolFont{cmcal}{OMS}{cmsy}{m}{n}
\SetSymbolFont{cmcal}{bold}{OMS}{cmsy}{b}{n}
\DeclareSymbolFontAlphabet{\mathcal}{cmcal}

\newcommand{\Z}[0]{\mathbb{Z}}

\title[]{A hypergeometric proof for a binomial identity related to $1/\pi$}

\author[B.~Hackl]{Benjamin Hackl}
\author[H.~Prodinger]{Helmut Prodinger}

\address[Benjamin Hackl]
{Institut f\"ur Mathematik,
  Alpen-Adria-Uni\-ver\-si\-t\"at Klagenfurt, Universit\"atsstra\ss e
  65--67, 9020 Klagenfurt, Austria}
\email{\href{mailto:benjamin.hackl@aau.at}{benjamin.hackl@aau.at}}
\thanks{B.~Hackl is supported by the Austrian
  Science Fund (FWF): P~28466-N35 and by a
  distinction grant from the Austrian Federal Ministry 
  of Education, Science and Research.}

\address[Helmut Prodinger]{Department of Mathematical
  Sciences, Stellenbosch University, 7602 Stellenbosch,
 South Africa}
\email{\href{mailto:hproding@sun.ac.za}{hproding@sun.ac.za}}


\begin{document}
\begin{abstract}
  We show that a binomial identity arising in the context of the study
  of series expansions of $1/\pi$ can be seen as an incarnation of
  Whipples second theorem for hypergeometric series.
\end{abstract}

\maketitle

In \cite{Sesma:2019:series}, Sesma studies a family of
series expansions for $1/\pi$ based on Heaviside's exponential
series. As a side product of these investigations, the hypergeometric
identity
\begin{equation}\label{eq:hyp}
  \sum_{n=0}^{m} \binom{m}{n} \binom{k+n}{m}
  \binom{k+m+n}{m+n}^{-1} \frac{2n + k + 1}{m + n + k + 1} = 1
\end{equation}
with $k \geq m \in \Z_{\geq 0}$ arises and was conjectured to be new.

We observe that rewriting the left-hand side of \eqref{eq:hyp} with
the help of factorials, regrouping the terms and simplifying quotients
with the help of the Pochhammer symbol $(\ell)_{n} = \ell (\ell +
1)\cdots (\ell + n - 1)$ yields
\begin{align}\notag
  \sum_{n=0}^{m}
  &\binom{m}{n} \binom{n+k}{m}
  \binom{m+n+k}{m+n}^{-1} \frac{2n + k + 1}{m + n + k + 1}\\ \notag
  &= \sum_{n=0}^{m} \frac{m!}{n!\,(m-n)!}
    \frac{(n+k)!}{m!\,(k-m+n)} \frac{(m+n)!\, k!}{(k + m + n)!}
    \frac{2n+k+1}{m+n+k+1}\\ \notag
  & = \frac{(k+1)!\, k!}{(k-m)!\, (k+m+1)!}
    \sum_{n=0}^{m} \frac{(k+n)!}{k!} \frac{2n+k+1}{k+1}
    \frac{(m+n)!}{(m-n)!} \frac{(k+m+n+1)!}{(k+m+n+1)!}
    \frac{(k-m)!}{(k-m+n)!} \frac{1}{n!}\\ \notag
  & = \frac{(k+1)!\, k!}{(k-m)!\, (k+m+1)!}
    \sum_{n\geq 0} (k+1)_{n} \frac{(1 +
    \frac{k+1}{2})_{n}}{(\frac{k+1}{2})_{n}}
    \frac{(-m)_{n} (-1)^{n} (m+1)_{n}}{(k+m+2)_{n} (k-m+1)_{n}}
    \frac{1}{n!}\\ \label{eq:simplification}
  & = \frac{(k+1)!\, k!}{(k-m)!\, (k+m+1)!}
    \pFq{4}{3}{k+1, 1 + \frac{k+1}{2}, -m, m+1}{\frac{k+1}{2}, k+m+2, k-m+1}{-1}.
\end{align}

The second Whipple theorem (cf.~\cite[Section 4.4]{Bailey:1964:hypergeometric}) for
$_{4}F_{3}$-hypergeometric series states that for suitable
complex-valued variables $a$, $b$, $c$ we have
\begin{equation}\label{eq:2nd-whipple}
  \pFq{4}{3}{a, 1+a/2, b, c}{a/2, a-b+1, a-c+1}{-1} =
  \frac{\Gamma(a-b+1) \Gamma(a-c+1)}{\Gamma(a+1)\Gamma(a-b-c+1)}.
\end{equation}
By comparing this to the hypergeometric $_{4}F_{3}$-term
in~\eqref{eq:simplification} we can immediately see that
we can use the second Whipple theorem with $a = k+1$, $b = -m$,
$c = m+1$ to rewrite~\eqref{eq:simplification} further to
\begin{align*}
  \frac{(k+1)!\, k!}{(k-m)!\, (k+m+1)!}
  \frac{\Gamma(k+m+2)\Gamma(k-m+1)}{\Gamma(k+2) \Gamma(k+1)} &=
  \frac{(k+1)!\, m!\, k!}{m!\, (k-m)!\, (m+k+1)!} \frac{(k+m+1)!\,
  (k-m)!}{(k+1)!\, k!}\\
  &=1,
\end{align*}
where we only had to use the relation $\Gamma(\ell+1) = \ell!$ for
non-negative integers $\ell$. This proves that the binomial
identity~\eqref{eq:hyp} is a consequence of the second
Whipple theorem.

\bibliographystyle{amsplainurl}
\bibliography{cheub,remark-1907.03188}

\end{document}